\documentclass[11pt,a4paper,reqno]{amsart}
\usepackage{amsmath}
\usepackage{amsfonts}
\usepackage{amssymb}
\usepackage{latexsym}

\numberwithin{equation}{section}
\newtheorem{proposition}{Proposition}[section]
\newtheorem{corollary}[proposition]{Corollary}
\newtheorem{theorem}[proposition]{Theorem}
\newtheorem{conjecture}[proposition]{Conjecture}

\newcommand{\In}{I_n}
\newcommand{\Jn}{J_n}
\newcommand{\Des}{Des}
\newcommand{\des}{d}
\newcommand{\fix}{fix}
\newcommand{\trans}{trans}
\newcommand{\maj}{maj}
\newcommand{\exc}{exc}

\newcommand{\wexc}{wexc}
\newcommand{\inv}{inv}
\newcommand{\DP}[2]{{\mathcal{I}}_{#1}^{\des}({#2})}
\newcommand{\SDP}[2]{\mathcal{J}_{#1}^{\des}(#2)}
\newcommand{\MP}[2]{{\mathcal{I}}_{#1}^{\maj}({#2})}
\newcommand{\EP}[2]{{\mathcal{I}}_{#1}^{\exc}({#2})}
\newcommand{\WEP}[2]{{\mathcal{I}}_{#1}^{\wexc}({#2})}
\newcommand{\IP}[2]{{\mathcal{I}}_{#1}^{\inv}({#2})}

\newcommand{\moloko}[3]{\overline{#1}^{(#2,#3)}} 
\newcommand{\molokt}[5]{\overline{#1}^{(#2,#3)(#4,#5)}}
\newcommand{\kodo}[3]{\moloko{#1}{#2}{#3}}
\newcommand{\kodt}[3]{\moloko{#1}{#2}{#3}}
\newcommand{\POLY}[3]{\mathcal{I}_{#2}^{#1}({#3})}
\newcommand{\JPOLY}[3]{\mathcal{J}_{#2}^{#1}({#3})}
\newcommand{\JPINV}[2]{\JPOLY{\inv}{#1}{#2}}
\newcommand{\sgn}{sgn}
\newcommand{\abs}[1]{|#1|}
\newcommand{\indicator}[1]{\mathbf{1}[#1]}

\newcommand{\TDIPB}[2]{\mathcal{IB}_{#1}(#2)}
\newcommand{\TDIPD}[2]{\mathcal{ID}_{#1}(#2)}
\newcommand{\TDIPO}[2]{\mathcal{IO}_{#1}(#2)}

\newcommand{\binv}{\inv_B}
\newcommand{\Jaut}{\mathcal{Y}}
\newcommand{\Amatrix}{\mathbf{A}}
\newcommand{\Wmatrix}{\mathbf{V}}
\newcommand{\DES}{D\'{e}sarm\'{e}nien}
\newcommand{\thediagram}{
\setlength{\unitlength}{0.240900pt}
\ifx\plotpoint\undefined\newsavebox{\plotpoint}\fi
\begin{picture}(1500,900)(0,0)
\font\gnuplot=cmr10 at 10pt
\gnuplot
\sbox{\plotpoint}{\rule[-0.200pt]{0.400pt}{0.400pt}}
\put(140.0,82.0){\rule[-0.200pt]{4.818pt}{0.400pt}}
\put(120,82){\makebox(0,0)[r]{ 0}}
\put(1419.0,82.0){\rule[-0.200pt]{4.818pt}{0.400pt}}
\put(140.0,212.0){\rule[-0.200pt]{4.818pt}{0.400pt}}
\put(120,212){\makebox(0,0)[r]{ 100}}
\put(1419.0,212.0){\rule[-0.200pt]{4.818pt}{0.400pt}}
\put(140.0,341.0){\rule[-0.200pt]{4.818pt}{0.400pt}}
\put(120,341){\makebox(0,0)[r]{ 200}}
\put(1419.0,341.0){\rule[-0.200pt]{4.818pt}{0.400pt}}
\put(140.0,471.0){\rule[-0.200pt]{4.818pt}{0.400pt}}
\put(120,471){\makebox(0,0)[r]{ 300}}
\put(1419.0,471.0){\rule[-0.200pt]{4.818pt}{0.400pt}}
\put(140.0,601.0){\rule[-0.200pt]{4.818pt}{0.400pt}}
\put(120,601){\makebox(0,0)[r]{ 400}}
\put(1419.0,601.0){\rule[-0.200pt]{4.818pt}{0.400pt}}
\put(140.0,730.0){\rule[-0.200pt]{4.818pt}{0.400pt}}
\put(120,730){\makebox(0,0)[r]{ 500}}
\put(1419.0,730.0){\rule[-0.200pt]{4.818pt}{0.400pt}}
\put(140.0,860.0){\rule[-0.200pt]{4.818pt}{0.400pt}}
\put(120,860){\makebox(0,0)[r]{ 600}}
\put(1419.0,860.0){\rule[-0.200pt]{4.818pt}{0.400pt}}
\put(140.0,82.0){\rule[-0.200pt]{0.400pt}{4.818pt}}
\put(140,41){\makebox(0,0){ 0}}
\put(140.0,840.0){\rule[-0.200pt]{0.400pt}{4.818pt}}
\put(284.0,82.0){\rule[-0.200pt]{0.400pt}{4.818pt}}
\put(284,41){\makebox(0,0){ 5}}
\put(284.0,840.0){\rule[-0.200pt]{0.400pt}{4.818pt}}
\put(429.0,82.0){\rule[-0.200pt]{0.400pt}{4.818pt}}
\put(429,41){\makebox(0,0){ 10}}
\put(429.0,840.0){\rule[-0.200pt]{0.400pt}{4.818pt}}
\put(573.0,82.0){\rule[-0.200pt]{0.400pt}{4.818pt}}
\put(573,41){\makebox(0,0){ 15}}
\put(573.0,840.0){\rule[-0.200pt]{0.400pt}{4.818pt}}
\put(717.0,82.0){\rule[-0.200pt]{0.400pt}{4.818pt}}
\put(717,41){\makebox(0,0){ 20}}
\put(717.0,840.0){\rule[-0.200pt]{0.400pt}{4.818pt}}
\put(862.0,82.0){\rule[-0.200pt]{0.400pt}{4.818pt}}
\put(862,41){\makebox(0,0){ 25}}
\put(862.0,840.0){\rule[-0.200pt]{0.400pt}{4.818pt}}
\put(1006.0,82.0){\rule[-0.200pt]{0.400pt}{4.818pt}}
\put(1006,41){\makebox(0,0){ 30}}
\put(1006.0,840.0){\rule[-0.200pt]{0.400pt}{4.818pt}}
\put(1150.0,82.0){\rule[-0.200pt]{0.400pt}{4.818pt}}
\put(1150,41){\makebox(0,0){ 35}}
\put(1150.0,840.0){\rule[-0.200pt]{0.400pt}{4.818pt}}
\put(1295.0,82.0){\rule[-0.200pt]{0.400pt}{4.818pt}}
\put(1295,41){\makebox(0,0){ 40}}
\put(1295.0,840.0){\rule[-0.200pt]{0.400pt}{4.818pt}}
\put(1439.0,82.0){\rule[-0.200pt]{0.400pt}{4.818pt}}
\put(1439,41){\makebox(0,0){ 45}}
\put(1439.0,840.0){\rule[-0.200pt]{0.400pt}{4.818pt}}
\put(140.0,82.0){\rule[-0.200pt]{312.929pt}{0.400pt}}
\put(1439.0,82.0){\rule[-0.200pt]{0.400pt}{187.420pt}}
\put(140.0,860.0){\rule[-0.200pt]{312.929pt}{0.400pt}}
\put(140.0,82.0){\rule[-0.200pt]{0.400pt}{187.420pt}}
\put(1259,815){\makebox(0,0)[r]{odd coefficients }}
\put(1279.0,815.0){\rule[-0.200pt]{28.908pt}{0.400pt}}
\put(169,94){\usebox{\plotpoint}}
\multiput(169.00,94.58)(0.659,0.498){85}{\rule{0.627pt}{0.120pt}}
\multiput(169.00,93.17)(56.698,44.000){2}{\rule{0.314pt}{0.400pt}}
\multiput(227.58,138.00)(0.499,0.623){111}{\rule{0.120pt}{0.598pt}}
\multiput(226.17,138.00)(57.000,69.758){2}{\rule{0.400pt}{0.299pt}}
\multiput(284.58,209.00)(0.499,0.751){113}{\rule{0.120pt}{0.700pt}}
\multiput(283.17,209.00)(58.000,85.547){2}{\rule{0.400pt}{0.350pt}}
\multiput(342.58,296.00)(0.499,0.768){113}{\rule{0.120pt}{0.714pt}}
\multiput(341.17,296.00)(58.000,87.518){2}{\rule{0.400pt}{0.357pt}}
\multiput(400.58,385.00)(0.499,0.733){113}{\rule{0.120pt}{0.686pt}}
\multiput(399.17,385.00)(58.000,83.576){2}{\rule{0.400pt}{0.343pt}}
\multiput(458.58,470.00)(0.499,0.614){111}{\rule{0.120pt}{0.591pt}}
\multiput(457.17,470.00)(57.000,68.773){2}{\rule{0.400pt}{0.296pt}}
\multiput(515.00,540.58)(0.580,0.498){97}{\rule{0.564pt}{0.120pt}}
\multiput(515.00,539.17)(56.829,50.000){2}{\rule{0.282pt}{0.400pt}}
\multiput(573.00,590.58)(1.168,0.497){47}{\rule{1.028pt}{0.120pt}}
\multiput(573.00,589.17)(55.866,25.000){2}{\rule{0.514pt}{0.400pt}}
\multiput(631.00,613.95)(12.518,-0.447){3}{\rule{7.700pt}{0.108pt}}
\multiput(631.00,614.17)(41.018,-3.000){2}{\rule{3.850pt}{0.400pt}}
\multiput(688.00,610.92)(0.785,-0.498){71}{\rule{0.727pt}{0.120pt}}
\multiput(688.00,611.17)(56.491,-37.000){2}{\rule{0.364pt}{0.400pt}}
\multiput(746.58,572.90)(0.499,-0.508){113}{\rule{0.120pt}{0.507pt}}
\multiput(745.17,573.95)(58.000,-57.948){2}{\rule{0.400pt}{0.253pt}}
\multiput(804.58,513.58)(0.499,-0.603){113}{\rule{0.120pt}{0.583pt}}
\multiput(803.17,514.79)(58.000,-68.790){2}{\rule{0.400pt}{0.291pt}}
\multiput(862.58,443.49)(0.499,-0.632){111}{\rule{0.120pt}{0.605pt}}
\multiput(861.17,444.74)(57.000,-70.744){2}{\rule{0.400pt}{0.303pt}}
\multiput(919.58,371.58)(0.499,-0.603){113}{\rule{0.120pt}{0.583pt}}
\multiput(918.17,372.79)(58.000,-68.790){2}{\rule{0.400pt}{0.291pt}}
\multiput(977.58,301.84)(0.499,-0.525){113}{\rule{0.120pt}{0.521pt}}
\multiput(976.17,302.92)(58.000,-59.919){2}{\rule{0.400pt}{0.260pt}}
\multiput(1035.00,241.92)(0.568,-0.498){99}{\rule{0.555pt}{0.120pt}}
\multiput(1035.00,242.17)(56.848,-51.000){2}{\rule{0.277pt}{0.400pt}}
\multiput(1093.00,190.92)(0.772,-0.498){71}{\rule{0.716pt}{0.120pt}}
\multiput(1093.00,191.17)(55.513,-37.000){2}{\rule{0.358pt}{0.400pt}}
\multiput(1150.00,153.92)(1.041,-0.497){53}{\rule{0.929pt}{0.120pt}}
\multiput(1150.00,154.17)(56.073,-28.000){2}{\rule{0.464pt}{0.400pt}}
\multiput(1208.00,125.92)(1.545,-0.495){35}{\rule{1.321pt}{0.119pt}}
\multiput(1208.00,126.17)(55.258,-19.000){2}{\rule{0.661pt}{0.400pt}}
\multiput(1266.00,106.92)(2.281,-0.493){23}{\rule{1.885pt}{0.119pt}}
\multiput(1266.00,107.17)(54.088,-13.000){2}{\rule{0.942pt}{0.400pt}}
\multiput(1324.00,93.93)(3.730,-0.488){13}{\rule{2.950pt}{0.117pt}}
\multiput(1324.00,94.17)(50.877,-8.000){2}{\rule{1.475pt}{0.400pt}}
\multiput(1381.00,85.94)(8.377,-0.468){5}{\rule{5.900pt}{0.113pt}}
\multiput(1381.00,86.17)(45.754,-4.000){2}{\rule{2.950pt}{0.400pt}}
\put(169,94){\raisebox{-.8pt}{\makebox(0,0){$\Diamond$}}}
\put(227,138){\raisebox{-.8pt}{\makebox(0,0){$\Diamond$}}}
\put(284,209){\raisebox{-.8pt}{\makebox(0,0){$\Diamond$}}}
\put(342,296){\raisebox{-.8pt}{\makebox(0,0){$\Diamond$}}}
\put(400,385){\raisebox{-.8pt}{\makebox(0,0){$\Diamond$}}}
\put(458,470){\raisebox{-.8pt}{\makebox(0,0){$\Diamond$}}}
\put(515,540){\raisebox{-.8pt}{\makebox(0,0){$\Diamond$}}}
\put(573,590){\raisebox{-.8pt}{\makebox(0,0){$\Diamond$}}}
\put(631,615){\raisebox{-.8pt}{\makebox(0,0){$\Diamond$}}}
\put(688,612){\raisebox{-.8pt}{\makebox(0,0){$\Diamond$}}}
\put(746,575){\raisebox{-.8pt}{\makebox(0,0){$\Diamond$}}}
\put(804,516){\raisebox{-.8pt}{\makebox(0,0){$\Diamond$}}}
\put(862,446){\raisebox{-.8pt}{\makebox(0,0){$\Diamond$}}}
\put(919,374){\raisebox{-.8pt}{\makebox(0,0){$\Diamond$}}}
\put(977,304){\raisebox{-.8pt}{\makebox(0,0){$\Diamond$}}}
\put(1035,243){\raisebox{-.8pt}{\makebox(0,0){$\Diamond$}}}
\put(1093,192){\raisebox{-.8pt}{\makebox(0,0){$\Diamond$}}}
\put(1150,155){\raisebox{-.8pt}{\makebox(0,0){$\Diamond$}}}
\put(1208,127){\raisebox{-.8pt}{\makebox(0,0){$\Diamond$}}}
\put(1266,108){\raisebox{-.8pt}{\makebox(0,0){$\Diamond$}}}
\put(1324,95){\raisebox{-.8pt}{\makebox(0,0){$\Diamond$}}}
\put(1381,87){\raisebox{-.8pt}{\makebox(0,0){$\Diamond$}}}
\put(1439,83){\raisebox{-.8pt}{\makebox(0,0){$\Diamond$}}}
\put(1339,815){\raisebox{-.8pt}{\makebox(0,0){$\Diamond$}}}
\put(1259,765){\makebox(0,0)[r]{even coefficients}}
\put(1279.0,765.0){\rule[-0.200pt]{28.908pt}{0.400pt}}
\put(140,83){\usebox{\plotpoint}}
\multiput(140.00,83.58)(0.831,0.498){67}{\rule{0.763pt}{0.120pt}}
\multiput(140.00,82.17)(56.417,35.000){2}{\rule{0.381pt}{0.400pt}}
\multiput(198.00,118.58)(0.606,0.498){91}{\rule{0.585pt}{0.120pt}}
\multiput(198.00,117.17)(55.786,47.000){2}{\rule{0.293pt}{0.400pt}}
\multiput(255.58,165.00)(0.499,0.560){113}{\rule{0.120pt}{0.548pt}}
\multiput(254.17,165.00)(58.000,63.862){2}{\rule{0.400pt}{0.274pt}}
\multiput(313.58,230.00)(0.499,0.725){113}{\rule{0.120pt}{0.679pt}}
\multiput(312.17,230.00)(58.000,82.590){2}{\rule{0.400pt}{0.340pt}}
\multiput(371.58,314.00)(0.499,0.838){113}{\rule{0.120pt}{0.769pt}}
\multiput(370.17,314.00)(58.000,95.404){2}{\rule{0.400pt}{0.384pt}}
\multiput(429.58,411.00)(0.499,0.905){111}{\rule{0.120pt}{0.823pt}}
\multiput(428.17,411.00)(57.000,101.292){2}{\rule{0.400pt}{0.411pt}}
\multiput(486.58,514.00)(0.499,0.838){113}{\rule{0.120pt}{0.769pt}}
\multiput(485.17,514.00)(58.000,95.404){2}{\rule{0.400pt}{0.384pt}}
\multiput(544.58,611.00)(0.499,0.707){113}{\rule{0.120pt}{0.666pt}}
\multiput(543.17,611.00)(58.000,80.619){2}{\rule{0.400pt}{0.333pt}}
\multiput(602.00,693.58)(0.537,0.498){105}{\rule{0.530pt}{0.120pt}}
\multiput(602.00,692.17)(56.901,54.000){2}{\rule{0.265pt}{0.400pt}}
\multiput(660.00,747.58)(1.370,0.496){39}{\rule{1.186pt}{0.119pt}}
\multiput(660.00,746.17)(54.539,21.000){2}{\rule{0.593pt}{0.400pt}}
\multiput(717.00,766.92)(1.842,-0.494){29}{\rule{1.550pt}{0.119pt}}
\multiput(717.00,767.17)(54.783,-16.000){2}{\rule{0.775pt}{0.400pt}}
\multiput(775.00,750.92)(0.580,-0.498){97}{\rule{0.564pt}{0.120pt}}
\multiput(775.00,751.17)(56.829,-50.000){2}{\rule{0.282pt}{0.400pt}}
\multiput(833.58,699.32)(0.499,-0.681){113}{\rule{0.120pt}{0.645pt}}
\multiput(832.17,700.66)(58.000,-77.662){2}{\rule{0.400pt}{0.322pt}}
\multiput(891.58,619.79)(0.499,-0.844){111}{\rule{0.120pt}{0.774pt}}
\multiput(890.17,621.39)(57.000,-94.394){2}{\rule{0.400pt}{0.387pt}}
\multiput(948.58,523.64)(0.499,-0.890){113}{\rule{0.120pt}{0.810pt}}
\multiput(947.17,525.32)(58.000,-101.318){2}{\rule{0.400pt}{0.405pt}}
\multiput(1006.58,420.81)(0.499,-0.838){113}{\rule{0.120pt}{0.769pt}}
\multiput(1005.17,422.40)(58.000,-95.404){2}{\rule{0.400pt}{0.384pt}}
\multiput(1064.58,324.08)(0.499,-0.755){111}{\rule{0.120pt}{0.704pt}}
\multiput(1063.17,325.54)(57.000,-84.540){2}{\rule{0.400pt}{0.352pt}}
\multiput(1121.58,238.70)(0.499,-0.569){113}{\rule{0.120pt}{0.555pt}}
\multiput(1120.17,239.85)(58.000,-64.848){2}{\rule{0.400pt}{0.278pt}}
\multiput(1179.00,173.92)(0.645,-0.498){87}{\rule{0.616pt}{0.120pt}}
\multiput(1179.00,174.17)(56.722,-45.000){2}{\rule{0.308pt}{0.400pt}}
\multiput(1237.00,128.92)(1.080,-0.497){51}{\rule{0.959pt}{0.120pt}}
\multiput(1237.00,129.17)(56.009,-27.000){2}{\rule{0.480pt}{0.400pt}}
\multiput(1295.00,101.92)(1.934,-0.494){27}{\rule{1.620pt}{0.119pt}}
\multiput(1295.00,102.17)(53.638,-15.000){2}{\rule{0.810pt}{0.400pt}}
\multiput(1352.00,86.93)(6.387,-0.477){7}{\rule{4.740pt}{0.115pt}}
\multiput(1352.00,87.17)(48.162,-5.000){2}{\rule{2.370pt}{0.400pt}}
\put(140,83){\makebox(0,0){$+$}}
\put(198,118){\makebox(0,0){$+$}}
\put(255,165){\makebox(0,0){$+$}}
\put(313,230){\makebox(0,0){$+$}}
\put(371,314){\makebox(0,0){$+$}}
\put(429,411){\makebox(0,0){$+$}}
\put(486,514){\makebox(0,0){$+$}}
\put(544,611){\makebox(0,0){$+$}}
\put(602,693){\makebox(0,0){$+$}}
\put(660,747){\makebox(0,0){$+$}}
\put(717,768){\makebox(0,0){$+$}}
\put(775,752){\makebox(0,0){$+$}}
\put(833,702){\makebox(0,0){$+$}}
\put(891,623){\makebox(0,0){$+$}}
\put(948,527){\makebox(0,0){$+$}}
\put(1006,424){\makebox(0,0){$+$}}
\put(1064,327){\makebox(0,0){$+$}}
\put(1121,241){\makebox(0,0){$+$}}
\put(1179,175){\makebox(0,0){$+$}}
\put(1237,130){\makebox(0,0){$+$}}
\put(1295,103){\makebox(0,0){$+$}}
\put(1352,88){\makebox(0,0){$+$}}
\put(1410,83){\makebox(0,0){$+$}}
\put(1339,765){\makebox(0,0){$+$}}
\end{picture}
	}

\title{Permutation statistics on involutions}
\author{W. M. B. Dukes}
\email{dukes@labri.fr}
\thanks{Supported by EC's Research Training Network `Algebraic Combinatorics in Europe', 
grant HPRN-CT-2001-00272 while the author was at Universit\`{a} di Roma Tor Vergata, Italy.}
\address{LaBRI, Universit\'{e} Bordeaux 1, 351 cours de la Lib\'{e}ration, 33405 Talence Cedex, France.}
\begin{document}
\maketitle
\begin{abstract}
In this paper we look at polynomials arising from statistics on the classes of involutions,
$I_n$, and involutions with no fixed points, $J_n$, in the symmetric group.
Our results are motivated by F. Brenti's conjecture~\cite{brenticonj} which states that the Eulerian distribution of $I_n$ is log-concave.
Symmetry of the generating functions is shown for the statistics $\des,\maj$ and the joint distribution $(\des,\maj)$.
We show that $\exc$ is log-concave on $I_n$, $\inv$ is log-concave on $J_n$ and $\des$ is partially unimodal on both $I_n$ and $J_n$.
We also give recurrences and explicit forms for the generating functions of the inversions statistic on involutions 
in Coxeter groups of types $B_n$ and $D_n$. Symmetry and unimodality of $inv$ is shown on the subclass 
of signed permutations in $D_n$ with no fixed points. 
In light of these new results, we present further conjectures at the end of the paper.
\end{abstract}

\section{Introduction}
In this paper we look at polynomials arising from statistics on the classes of involutions and involutions with no fixed points
in the symmetric group. 

Let $S_n$ be the symmetric group on $[1,n]$.
Call $\Des(\sigma) := \{\,i\,:\,1\leq i <n \mbox{ and } \sigma_i > \sigma_{i+1}\}$
the {\it{descent set}} of $\sigma \in S_n$ and 
the number of {\it{descents}} is denoted $\des(\sigma):=|\Des(\sigma)|$.
We further define $\des_i(\sigma):= |\{j\geq i:j\in\Des(\sigma)\}|$, the {\it{partial descents}} of $\sigma$ for $1\leq i<n$.
The {\it{major index}} of $\sigma$ is $\maj(\sigma) := \sum_{i\in \Des(\sigma)} i$
and the number of {\it{inversions}} is $inv(\sigma):= |\{1\leq i<j\leq n: \sigma_i>\sigma_j\}|$.
The number of {\it{excedances}} is $\exc(\sigma) := |\{1\leq i\leq n: \sigma_i>i\}|$  and 
{\it{weak excedances}} is $\wexc(\sigma):= |\{1\leq i\leq n: \sigma_i\geq i\}|$.
Let $\fix(\sigma)$ and $\trans(\sigma)$ denote the number of fixed points and transpositions of $\sigma$, respectively.
We use the notation $[x^i]P(x)$ for the coefficient of $x^i$ in the polynomial $P(x)$.

For a statistic $stat: S_n \to {\mathbf{N}}_0$, define the polynomials
\begin{eqnarray*}
  \POLY{stat}{n}{x} \;:=\; \sum_{\sigma \in \In} x^{stat(\sigma)}, &&
  \JPOLY{stat}{n}{x} \;:=\; \sum_{\sigma \in \Jn} x^{stat(\sigma)},
\end{eqnarray*}
where $\In := \{\sigma\in S_n : \sigma^2 = \mbox{id}\}$ and $\Jn:= \{\sigma \in \In : \fix(\sigma)=0\}$. 
For an arbitrary collection $S_n' \subseteq S_n$, the sequence of coefficients of  
$\sum_{\pi \in S_n'} x^{\des(\pi)}$ is termed the {\it{Eulerian distribution}} of $S_n'$.
The results in this paper are motivated by 
\begin{conjecture}[Brenti~\cite{brenticonj}]
The Eulerian distribution of $I_n$ is log-concave.
\end{conjecture}
We propose further conjectures concerning statistics on different classes of involutions in
the final section.

\section{Involutions in the Symmetric group}
\subsection{The excedances statistic}
\begin{theorem}
The coefficients of the polynomial $\EP{n}{x}$ are log-concave.
\end{theorem}

\begin{proof}
The number $\exc(\sigma)$ is precisely the number of 2-cycles in an involution, so we have 
\begin{eqnarray}
\label{exceqn}
\EP{n}{x} &=& \sum_{k=0}^{\lfloor n/2 \rfloor} \dfrac{n!}{k! (n-2k)!} \left( \dfrac{x}{2}\right)^k.
\end{eqnarray}
It is an easy exercise to show log-concavity for $0\leq j < \lfloor n/2 \rfloor$ since we have a direct expression for the coefficients.
\end{proof}
Note that the polynomials $\EP{n}{x}$ are closely related to the Hermite polynomials $h_n(x)$, whereby
\begin{eqnarray*}
\sum_{n\geq 0 } \dfrac{h_n(x) t^n}{n!} &=& \exp ( tx-t^2/2),
\end{eqnarray*}
via the equation $\EP{n}{x} = (-x)^n h_n(-1/2x)$. 
The Hermite polynomials are known to be real-rooted (see for example Stanley~\cite[p. 505]{stanley.ann}).

The Sch\"{u}tzenberger involution on tableaux, $T \to evac (T)$, maps involutions to involutions and
$\wexc(evac(\sigma)) = n-\exc(\sigma)$, since $evac(\sigma)_i = n+1-\sigma_{n+1-i}$,
so that
$\WEP{n}{x} = x^{n} \EP{n}{x^{-1}}$, hence

\begin{corollary}
The coefficients of the polynomial $\WEP{n}{x}$ are log-concave.
\end{corollary}

\subsection{The descents and major index statistics}
In the spirit of Adin et. al.~\cite{abr}, we define 
\begin{eqnarray*}
\mathcal{G}_n(x_1,\ldots ,x_{n-1}) & := & 
	\sum_{\sigma \in\In} x_1^{d_1(\sigma)} x_2^{d_2(\sigma)} \cdots x_{n-1}^{d_{n-1}(\sigma)}.
\end{eqnarray*}

\begin{theorem}
\label{g-thm}
The polynomial $\mathcal{G}_n(x_1,\ldots ,x_{n-1})$ satisfies
\begin{eqnarray*}
  \mathcal{G}_n(x_1,\ldots ,x_{n-1}) &=& x_1^{n-1}x_2^{n-2} \cdots x_{n-1} \mathcal{G}_n(x_1^{-1},\ldots ,x_{n-1}^{-1}).
\end{eqnarray*}
\end{theorem}

\begin{proof}
If $\sigma \in \In$ then the reading and insertion tableau associated with $\sigma$ under Robinson-Schensted 
correspondence (Stanley~\cite[Ch. 7]{ec2}) are identical. 
That is, there is a bijection between $\In$ and all standard Young tableaux (SYT) on $[1,n]$.

Let $\sigma \in \In$ with associated SYT $T$. 
The set $\Des(\sigma)$ corresponds to those entries $i$ in the tableau $T$ such that $(i+1)$ is below and weakly to the left of $i$. 
Let $T^{\perp}$ be the tableau $T$ reflected on its main diagonal. 
Notice that if $(i+1)$ is below and weakly to the left of $i$ in $T$, 
then $(i+1)$ is to the right of and weakly above $i$ in $T^{\perp}$. 
The bijection between the class of  SYT on $[1,n]$ and involutions $\In$ shows that
to $T^{\perp}$ there corresponds a unique involution $\sigma^{\perp} \in \In$, and has the property that 
$\{\Des(\sigma),\Des(\sigma^{\perp})\}$ is a partition of the set $[1,n-1]$.
In this manner, the reflection operation is an involution on involutions.

It follows that
\begin{eqnarray*}
  d_i(\sigma^{\perp}) &=& |\{j\geq i: j\in Des(\sigma^{\perp})\}|  \\
  &=& |\{j\geq i : j\not\in \Des(\sigma)\}| \\
  &=& n-i - |\{j\geq i : j\in \Des(\sigma)\}| \\
  &=& n-i - d_i(\sigma).
\end{eqnarray*}

We have shown that if $\sigma \in \In$, 
then there is a unique $\sigma^{\perp} \in \In$ such that 
$(d_1(\sigma^{\perp}),\ldots , d_{n-1}(\sigma^{\perp})) = (n-1-d_1(\sigma),\ldots , 1-d_{n-1}(\sigma))$.
\end{proof}

Both polynomials $\DP{n}{q}$ and $\MP{n}{q}$ are instances of the $\mathcal{G}$ polynomial since
$\DP{n}{q} = \mathcal{G}_n(q,1,\ldots ,1)$ and 
$\MP{n}{q} = \mathcal{G}_n(q,q,\ldots ,q)$. 
Comparing coefficients on both sides of the symmetric $\mathcal{G}$ relation yields

\begin{corollary}
The polynomials $\DP{n}{t}$ and $\MP{n}{t}$ are symmetric.
\end{corollary}

Symmetry of the polynomials $\DP{n}{x}$ and $\JPOLY{\des}{n}{t}$ was conjectured by Dumont and first proven by Strehl~\cite{strehl}, 
using a method similar to that of the previous theorem for the coefficients of $\DP{n}{x}$. 
A separate argument was used to prove symmetry of $\JPOLY{\des}{n}{t}$ 
because for $\sigma \in J_n$, it is not necessarily true that $\sigma^{\perp}\in J_n$.
Theorem~\ref{g-thm} allows us to
show symmetry of the joint distribution of $(\des,\maj)$ on $\In$ since
$\sum_{\sigma\in\In} t^{\des(\sigma)} q^{\maj(\sigma)} = \mathcal{G}_n(tq,q,\ldots ,q)$.

\begin{corollary}
The polynomial
\begin{eqnarray*}
  \mathcal{I}_n^{\des,\maj}(t,q)&=&\sum_{\sigma\in\In} t^{\des(\sigma)} q^{\maj(\sigma)} 
\end{eqnarray*}
is symmetric in the sense that $[t^i q^j] \mathcal{I}_n^{\des,\maj}(t,q) = [t^{n-1-i}q^{{n\choose 2}-j}] \mathcal{I}_n^{\des,\maj}(t,q)$.
\end{corollary}

Hultman~\cite{hultman} recently proved that for any finite Coxeter system $(W,S)$, the associated descent polynomial 
$\sum_w t^{\des_W(w)}$ is symmetric
where the sum ranges over all $w \in W$ with $w^2= \mbox{{id}}_w$.
{\DES} and Foata~\cite{df} use an elegant method involving Schur functions to derive the generating function
\begin{eqnarray}
  \label{thm_df}
  \sum_{n\geq 0} \dfrac{H_n(z_1,z_2,t,q) u^n}{(t;q)_n} &=& 
  \sum_{r\geq 0} t^r \dfrac{1}{(z_1u;q)_{r+1}} \prod_{0\leq i<j\leq r} \dfrac{1}{1-u^2 z_2 q^{i+j}}\quad
\end{eqnarray}
where $H_n(z_1,z_2,t,q):= \sum_{\sigma \in I_n} z_1^{fix(\sigma)} z_2^{trans(\sigma)} t^{d(\sigma)} q^{maj(\sigma)}$,
$(a;q)_0=1$ and $(a;q)_n = (1-a) (1-aq) \cdots (1-aq^{n-1})$. 
The generating functions for the polynomials $\DP{n}{t}$, $\MP{n}{q}$ are immediate from this:
\begin{eqnarray}
\sum_{n\geq 0} \dfrac{\DP{n}{t} u^n}{(1-t)^n} &=& \sum_{r\geq 0} t^r \left( \dfrac{1}{(1-u)^{r+1}(1-u^2)^{r(r+1)/2}}   \right); \label{dgf} \\
\sum_{n\geq 0} \dfrac{\SDP{n}{t} u^n}{(1-t)^n} &=& \sum_{r\geq 0} t^r \left( \dfrac{1}{(1-u^2)^{r(r+1)/2}}   \right); \label{nofixedgf} \\
\sum_{n\geq 0} \dfrac{\MP{n}{q} u^n}{(q;q)_{n}} &=& \sum_{r\geq 0}  \dfrac{1}{(u;q)_{r+1}} \prod_{0\leq i<j\leq r} \dfrac{1}{(1-u^2q^{i+j})}.
\end{eqnarray}
By extracting the appropriate coefficients, we now show partial unimodality of $\DP{n}{q}$ and $\SDP{n}{q}$.
The onerous aspect of proving total unimodality using these 
direct expressions seems to be the appearance of both $r$ and ${r+1\choose 2}$ within binomial terms.

\begin{theorem}
For all $1\leq i\leq n^{0.925}/10$, 
$[t^i]\JPOLY{\des}{n}{t}< [t^{i+1}]\JPOLY{\des}{n}{t}$ 
and
$[t^{n+1-i}]\JPOLY{\des}{n}{t}> [t^{n+2-i}]\JPOLY{\des}{n}{t}$.
\end{theorem}

\begin{proof}
Extracting the coefficient of $u^n$ in Equation (\ref{nofixedgf}), one finds 
\newcommand{\joc}[2]{\alpha_{#1,#2}}
\begin{eqnarray*}
\SDP{n}{t} &=& \sum_{p=1}^n \joc{n}{p}t^p \;=\;\sum_{p=1}^n t^p\left\{ \sum_{k=0}^{p-1} (-1)^k {n+1 \choose k} {{p-k+1 \choose 2}+n/2-1 \choose n/2}\right\}.
\end{eqnarray*}
Inverting this gives
\begin{eqnarray*}
f_n(p) := {{p+1\choose 2}+n/2-1\choose n/2} &=& \sum_{i=0}^{p-1} {n+i \choose n} \joc{n}{p-i}.
\end{eqnarray*}
For $p \geq 2$,
\begin{eqnarray*}
\lefteqn{f_n(p)-f_n(p-1) } \\
	&=& \joc{n}{p} - \joc{n}{p-1} 
	+ \sum_{i=1}^{p-1} {n+i\choose i} \joc{n}{p-i} - \sum_{i=1}^{p-2} {n+i\choose i} \joc{n}{p-1-i}\\
	&=& \joc{n}{p} - \joc{n}{p-1} 
	+{n+p-1\choose p-1}\joc{n}{1} + \sum_{i=1}^{p-2} {n+i\choose i} (\joc{n}{p-i}- \joc{n}{p-1-i})\\
	&\leq & \joc{n}{p} - \joc{n}{p-1} 
	+{n+p-1\choose p-1}\joc{n}{1} + (n+1)\sum_{i=1}^{p-2} {n+i\choose i} (\joc{n}{p-i}- \joc{n}{p-1-i})\\
	&\leq & \joc{n}{p} - \joc{n}{p-1} 
	+(n+1)(f_n(p-1)-f_n(p-2)).
\end{eqnarray*}
Thus $\joc{n}{p} - \joc{n}{p-1} \geq f_n(p)-f_n(p-1) - (n+1)(f_n(p-1)-f_n(p-2))$.
The right hand side of the previous inequality is positive for $p$ not too large. Notice that
\begin{eqnarray*}
\dfrac{f_n(p)}{f_n(p-1)} & \geq& \left(1+\dfrac{n}{p^2+p-2}\right)^p \\
&\geq & \left(1+\dfrac{n+2}{2p^2}\right)^p
\end{eqnarray*}
which, in turn, is bounded below by $n+2$ 
when $p\leq n^{0.925}/10$.
The second inequality follows from symmetry as shown in Strehl~\cite{strehl}.
\end{proof}

\begin{theorem}
For all $1\leq k \leq 0.175 n^{0.931}$,
$[t^{k-1}]\DP{n}{t}<[t^{k}]\DP{n}{t}$ and $[t^{n-1-k}]\DP{n}{t}>[t^{n-k}]\DP{n}{t}$.
\end{theorem}

\begin{proof}
Extracting the coefficient of $u^n$ in Equation (\ref{dgf}) we find
\begin{eqnarray}
  \DP{n}{t} &=& \sum_{k=0}^{n-1} \beta_{n,k} t^k \; = \; \sum_{k=0}^{n-1} t^k \left\{ 
                \sum_{j=0}^k {n+1 \choose j} (-1)^{j} \gamma (n,k-j)  \right\} \nonumber
\end{eqnarray}
where $\gamma(n,0);=1$ and $\gamma(n,r) := \sum_{i=0}^{\lfloor n/2\rfloor}  {n-2i+r \choose r} {i+{r+1\choose 2}-1 \choose i}$ for $r>0$.
Thus for $0\leq k \leq n$,
\begin{eqnarray*}
  \gamma(n,k) &=& \sum_{i=0}^{k} {n+i\choose i} \beta_{n,k-i}
\end{eqnarray*}
and so 
\begin{eqnarray*}
  \lefteqn{\gamma(n,k) - \gamma(n,k-1)}\\
  &=& \sum_{i=0}^{k-1} {n+i\choose i} (\beta_{n,k-i} -\beta_{n,k-1-i}) + {n+k\choose k} \\
  &=& \beta_{n,k} - \beta_{n,k-1} + \sum_{i=0}^{k-2} {n+i+1\choose i+1} (\beta_{n,k-1-i}-\beta_{n,k-2-i}) +{n+k\choose k}  \\
  &<& \beta_{n,k} - \beta_{n,k-1} + (n+1)\sum_{i=0}^{k-2} {n+i\choose i} (\beta_{n,k-1-i}-\beta_{n,k-2-i}) +{n+k\choose k} \\
  &=& \beta_{n,k} - \beta_{n,k-1} + (n+1) \left(\gamma(n,k-1)-\gamma(n,k-2)-{n+k-1\choose k-1}\right)\\ && \;+{n+k\choose k}  \\
  &<& \beta_{n,k} - \beta_{n,k-1} + (n+1) \gamma(n,k-1).
\end{eqnarray*}
It suffices to show that for $n$ and $k$ as stated in the theorem, $\gamma(n,k)> (n+2) \gamma(n,k-1)$.
One may also write 
\begin{eqnarray*}
  \gamma (n,k) &=& [u^n] \dfrac{(1+u)^{k+1}}{(1-u^2)^{k+2\choose 2}}
\end{eqnarray*}
and since $n >k+1$ we have
\begin{eqnarray*}
  \gamma(n,k) &=& \sum_{i=0}^{\lfloor k+1/2\rfloor} {k+1 \choose n-2\lfloor n/2\rfloor +2i} 
	[u^{2(\lfloor n/2\rfloor -i)}] \dfrac{1}{(1-u^2)^{k+2\choose 2}}\\
	&>& \sum_{i=0}^{\lfloor k/2\rfloor} {k \choose n-2\lfloor n/2\rfloor +2i} 
	[u^{2(\lfloor n/2\rfloor -i)}] \dfrac{1}{(1-u^2)^{k+2\choose 2}}.
\end{eqnarray*}
Now for all $m\geq (n-k-2)/2$, 
\begin{eqnarray*}
  {{k+2\choose 2} + m\choose m+1} & > & 
	\left(1+\dfrac{2m}{(k+1)(k+2)}  \right)^{k+1}  {{k+1\choose 2} + m\choose m+1} \\
	&>& (n+2) {{k+1\choose 2} + m\choose m+1}.
\end{eqnarray*}
for $k \leq 0.175 n^{0.931}$, hence
\begin{eqnarray*}
  \gamma(n,k) & > & (n+2) 
	\sum_{i=0}^{\lfloor k/2\rfloor} {k \choose n-2\lfloor n/2\rfloor +2i} 
	    [u^{2(\lfloor n/2\rfloor -i)}] \dfrac{1}{(1-u^2)^{k+1\choose 2}}\\
&=& (n+2) \gamma(n,k-1),
\end{eqnarray*}
giving the first inequality.
Again, symmetry of the $\DP{n}{t}$ polynomial yields the second inequality.
\end{proof}

\subsection{The inversions statistic}
The generating function for the inversions statistic on involutions 
is intimately related to the $q$-Hermite polynomials, as studied by \DES~\cite{jd}.
Let $a_{n}(k,j)$ be the number of involutions in $\In$ with $k$ fixed points and $j$ inversions,
and define $Z_n(x,q) := \sum_{k,j} a_{n}(k,j) q^j x^k$. 
\DES~\cite[Eqns. 3.10,3.11]{jd} showed 
\begin{eqnarray*}
Z_{n+1}(x,q) &=& xZ_n(x,q) + q\left(\dfrac{1-q^{2n}}{1-q^2}\right) Z_{n-1}(x,q) 
\end{eqnarray*}
for all $n>1$ with $Z_0(x,q)=1$ and $Z_1(x,q)=x$.
Setting $x=1,0$, yields the following proposition.

\begin{proposition}
\label{desprop}
For all $n\geq 0$,
\begin{eqnarray*} 
  \IP{n+2}{q} &=& \IP{n+1}{q} + q\left(\dfrac{1-q^{2(n+1)}}{1-q^2}  \right)\IP{n}{q}, 
\end{eqnarray*}
where $\IP{0}{q} , \IP{1}{q} :=1$ and
for $n\geq 0$,
\begin{eqnarray*}
  \JPINV{n+2}{q} &=& q\left(\dfrac{1-q^{2(n+1)}}{1-q^2}\right)\JPINV{n}{q},
\end{eqnarray*}
where $\JPINV{0}{q} = 1$.
\end{proposition}

The above recurrences can also be derived in a straightforward manner using a special case of Equation (\ref{bcon}).
The coefficients of $\IP{n}{q}$ are neither log-concave nor unimodal (see Figure 1) but the recursion in the previous proposition admits a solution as a matrix product, which may be of benefit in approaching Conjecture~\ref{myconjecture}(ii).

\begin{proposition}
\label{matrixrecur}
Let $g_i(q):=\sum_{j=0}^{i-2} q^{1+2j}$ and 
$\Amatrix_n(q):= \prod_{i=2}^n \left( \begin{array}{cc} 1 & 1 \\ g_{i}(q) & 0  \end{array}  \right)$,
then $\IP{n}{q} = \Amatrix_n(q)_{1,1} + \Amatrix_n(q)_{2,1}$.
\end{proposition}

\begin{proof}
For all $n\geq 2$, we may write
\begin{eqnarray*}
  \IP{n}{q} &=& \prod_{i=1}^n \Jaut_i(q)
\end{eqnarray*}
where $\Jaut_n(q):= \IP{n}{q} / \IP{n-1}{q}$ and $\IP{0}{q} , \IP{1}{q} :=1$. 
From the first recurrence in Proposition~\ref{desprop} the polynomial $\Jaut_n(q)$ satisfies the recurrence
$\Jaut_n(q) = 1+g_{n-1}(q)/\Jaut_{n-1}(q)$ for all $n\geq 2$ where $g_n(q):=q+q^3+\ldots + q^{2n-3}$.
Using this, the product $\Jaut_i \Jaut_{i+1} \cdots \Jaut_{n}$ may be written in the form $\alpha_i(q) \Jaut_i(q) + \beta_i(q)$.
It is easily seen that $\alpha_n(q) = 1$, $\beta_n(q) = 0$ and 
\begin{eqnarray*}
  \left(\begin{array}{c} \alpha_i(q) \\ \beta_i (q) \end{array}\right) 
  &=& \left( \begin{array}{cc} 1 & 1 \\ g_{i+1}(q) & 0  \end{array}  \right)
  \left(\begin{array}{c} \alpha_{i+1}(q) \\ \beta_{i+1}(q)\end{array}\right).
\end{eqnarray*}
Thus we have  $\IP{n}{q} = \alpha_1(q) \Jaut_1(q) + \beta_1(q) = \alpha_1(q) + \beta_1(q)$, since $\Jaut_1(q)=1$, and
\begin{eqnarray*}
  \left(\begin{array}{c} \alpha_1(q) \\ \beta_1(q)  \end{array}\right) &=&
  \left( \begin{array}{cc} 1 & 1 \\ g_{2}(q) & 0  \end{array}  \right)
  \left( \begin{array}{cc} 1 & 1 \\ g_{3}(q) & 0  \end{array}  \right) 
  \cdots
  \left( \begin{array}{cc} 1 & 1 \\ g_{n}(q) & 0  \end{array}  \right) 
  \left(\begin{array}{c} 1 \\ 0  \end{array}\right).
\end{eqnarray*}
\end{proof}

\begin{theorem}
The coefficients of the polynomial $\JPINV{n}{q} $ are log-concave.
\end{theorem}

\begin{proof}
Solving the second recurrence in Proposition~\ref{desprop} we get:
\begin{eqnarray*}
  \JPINV{2m}{q} &=& {q^m} \prod_{i=1}^{m-1} \dfrac{1-q^{2(2i+1)}}{1-q^2}.
\end{eqnarray*}
Set $u=q^2$ and notice that the sequence of non-zero coefficients in $\JPINV{n}{q}$ is
the same as $\prod_{i=1}^{m-1} \dfrac{1-u^{2i+1}}{1-u}$. 
The coefficients of the polynomials $(1-u^{2i+1})/(1-u)$ are 
non-negative log-concave sequences with no internal zero coefficients.
Thus using Stanley~\cite[Prop. 2]{stanley.ann}, the product of all 
such polynomials will also be log-concave with no internal zero coefficients.
\end{proof}

\section{Involutions in Coxeter groups of types $B$ and $D$}
In this section we give recursive expressions for the inversion polynomials of involutions 
for Coxeter groups of types $B$ and $D$. 
We use the notation of Bj\"{o}rner and Brenti~\cite{coxbook}.

Coxeter groups of type $B$, the `signed permutations', are defined as follows:
let $S_n^{B}$ be the group of all bijections $\pi$ on the set 
$[\pm n]\backslash \{0\}$ such that $\pi(-a)=-\pi(a)$ for all $a \in [\pm n]$.
For $\pi \in S_n^B$, define
\begin{eqnarray*}
  N_1 (\pi(1),\ldots ,\pi(n)) &:=& |\{1\leq i \leq n: \pi(i)<0\}| \\
  N_2 (\pi(1),\ldots ,\pi(n)) &:=& |\left\{ 1\leq i <j \leq n : \pi(i)+\pi(j) < 0 \right\}|.
\end{eqnarray*}
Let $S_n^D$ be the subgroup of $S_n^B$ consisting of all signed permutations $\pi \in S_n^B$ such that
there are an even number of negative entries in the window of $\pi$, i.e.
$S_n^D := \{ \pi \in S_n^B : N_1(\pi) \equiv 0 (\mbox{mod } 2)\}$.
For completeness let us also define those signed permutations containing an odd number of negative signs in 
the window of $\pi$, $S_n^O = S_n^{B} \backslash S_n^D$.

The inversions statistics on $S_n^B$ and $S_n^D$ are defined slightly differently to $inv$ on $S_n$. 
From \cite[Equations (8.1) and (8.18)]{coxbook}, 
\begin{eqnarray*}
  \binv (\pi) &:=& \inv (\pi(1),\ldots ,\pi(n)) + N_1 (\pi(1) ,\ldots ,\pi(n)) \\ 
              &  & + N_2 (\pi(1),\ldots ,\pi(n)) \\
  \inv_D(\pi) &:=& \inv(\pi(1),\ldots ,\pi(n)) + N_2 (\pi(1),\ldots ,\pi(n)).
\end{eqnarray*}
Let us mention that in the symmetric group setting, 
\begin{eqnarray*}
  \sum_{\pi \in S_n^B} q^{\inv_B(\pi)} &=& [2]_q [4]_q \ldots [2n]_q\\
  \sum_{\pi \in S_n^D} q^{\inv_D(\pi)} &=& [2]_q [4]_q \cdots [2n-2]_q [n]_q.
\end{eqnarray*}
where
$[i]_q := 1+q+q^2+ \ldots +q^{i-1}$ 
(see \cite[Theorem 7.1.5.]{coxbook})

Define $I_n^B := \{ \pi \in S_n^B : \pi^2 = \mbox{id} \}$,
$I_n^D : = \{ \pi \in S_n^D : \pi^2 = \mbox{id}\}$ and $I_n^{O}:=I_n^B\backslash I_n^D$.
Let 
\begin{eqnarray*}
\TDIPB{n}{q}  &:=& \sum_{\pi \in I_n^B} q^{\inv_B(\pi)},\\
\end{eqnarray*}
with $\TDIPD{n}{q}$ and $\TDIPO{n}{q}$ similarly defined.
To aid the proof of the following two theorems, we introduce some notation concerning 
the recursive construction of these signed permutations.

Let $\pi \in I_{n}^B$ and denote by $\kodo{\pi}{n+1}{n+1}$ the signed permutation $\pi ' \in I_{n+1}^B$ such that
$\pi'(i) = \pi(i)$, for $1\leq i, \leq n$ and $\pi'(n+1) = n+1$. 
Similarly let $\kodo{\pi}{-(n+1)}{n+1}$ be the signed permutation $\pi ' \in I_{n+1}^B$ such that
$\pi'(i) = \pi(i)$ for $1\leq i \leq n$ and $\pi'(n+1) = -(n+1)$.

For $\pi\in I_n^B$ and $k\in[\pm (n+1)]-\{0\}$, let $\kodt{\pi}{k}{n+2}$ be the signed permutation $\pi ' \in I_{n+2}^B$ such that
\begin{itemize}
\item $\pi '(\abs{k}) = (n+2) \sgn(k) $, $\pi ' (n+2) = k$,
\item for all $1 \leq i \leq n$, 
  \begin{eqnarray*}
    \pi '(i+\indicator{i\geq \abs{k}}) &=& \pi (i) + \sgn(\pi(i)) \indicator{\abs{\pi(i)} \geq \abs{k}}
  \end{eqnarray*}
\end{itemize}
where $\sgn(a)=+1$ if $a>0$ and $-1$ otherwise.
Consequently $I_{n+2}^B$, $I_{n+2}^D$ and $I_{n+2}^O$ may be constructed recursively,
\begin{eqnarray}
I_{n+2}^{B} &=& \biguplus_{\pi \in I_{n+1}^B} \{ \kodo{\pi}{n+2}{n+2} , \kodo{\pi}{-(n+2)}{n+2} \} \; \uplus \;  \nonumber\\ &&
	\biguplus_{k=1\atop \pi \in I_n^B}^{n+1} \{ \kodt{\pi}{k}{n+2} , \kodt{\pi}{-k}{n+2}   \}  \label{bcon}\\
I_{n+2}^{D} &=& 
	\biguplus_{\pi \in I_{n+1}^D} \{ \kodo{\pi}{n+2}{n+2} \} \; \uplus \; 
	\biguplus_{\pi \in I_{n+1}^O} \{ \kodo{\pi}{-(n+2)}{n+2} \} \; \uplus \;  \nonumber\\ &&
	\biguplus_{k=1\atop \pi \in I_n^D}^{n+1} \{ \kodt{\pi}{k}{n+2} , \kodt{\pi}{-k}{n+2}   \} \label{dcon} \\
I_{n+2}^{O} &=& 
	\biguplus_{\pi \in I_{n+1}^O} \{ \kodo{\pi}{n+2}{n+2} \} \; \uplus \; 
	\biguplus_{\pi \in I_{n+1}^D} \{ \kodo{\pi}{-(n+2)}{n+2} \} \; \uplus \;  \nonumber \\ &&
	\biguplus_{k=1\atop \pi \in I_n^O}^{n+1} \{ \kodt{\pi}{k}{n+2} , \kodt{\pi}{-k}{n+2}   \} . \label{ocon}
\end{eqnarray}

\begin{theorem} 
For all $n\geq 2$,
\begin{eqnarray*}
  \TDIPB{n+2}{q} &=& (1+q^{2n+3}) \TDIPB{n+1}{q} + \dfrac{q(1+q^2)(1-q^{2(n+1)})}{1-q^2} \TDIPB{n}{q}
\end{eqnarray*}
with initial polynomials 
$\TDIPB{2}{q} = 1+2q+2q^3+q^4$,
$\TDIPB{3}{q} = 1+3q+q^2+3q^3+2q^4+2q^5+3q^6+q^7+3q^8+q^9$.
\end{theorem}

\begin{proof}
Using Equation (\ref{bcon}), 
\begin{eqnarray*}
  {\TDIPB{n+2}{q}} &=& \sum_{\pi \in I_{n+1}^B} q^{inv_B(\kodo{\pi}{n+2}{n+2})} + q^{inv_B(\kodo{\pi}{-(n+2)}{n+2})} \\
  && + \sum_{k=1}^{n+1} \sum_{\pi \in I_n^B} q^{inv_B(\kodt{\pi}{k}{n+2})} + q^{inv_B(\kodt{\pi}{-k}{n+2})}.
\end{eqnarray*}

If $\pi \in I_{n+1}^B$, then $inv_B(\kodo{\pi}{n+2}{n+2}) = inv_B(\pi(1),\ldots ,\pi(n+1),n+2) = inv_B(\pi)$ 
and $inv_B(\kodo{\pi}{-(n+2)}{n+2}) = inv_B(\pi)+2n+3$.
Similarly if $\pi \in I_n^B$ and $1\leq k \leq n+1$, then $inv_B(\kodt{\pi}{k}{n+2}) = inv_B(\pi) +2n+3-2k$
and $inv_B(\kodt{\pi}{-k}{n+2}) = inv_B(\pi) + 2k+1$.  Hence
\begin{eqnarray*}
{\TDIPB{n+2}{q}} &=& \sum_{\pi \in I_{n+1}^B} q^{inv_B(\pi)} + q^{inv_B(\pi) +2n+3} \\
&& + \sum_{k=1}^{n+1} \sum_{\pi \in I_n^B} q^{inv_B(\pi)+2n-2k+3} + q^{inv_B(\pi)+2k+1}\\
&=& (1+q^{2n+3}) \TDIPB{n+1}{q} + \TDIPB{n}{q} \sum_{k=1}^{n+1} (q^{2n-2k+3}+q^{2k+1}).
\end{eqnarray*}
\end{proof}

We may express $\TDIPB{n}{q}$ in a somewhat closed form, as was done in {{Proposition~\ref{matrixrecur}}};
for all $n\geq 3$,
$\TDIPB{n}{q} = (\Wmatrix_n(q)_{1,1}+\Wmatrix_n(q)_{2,1})(1+2q+2q^3+q^4)$  where
\begin{eqnarray*}
\Wmatrix_n(q) &=& \prod_{i=3}^{n} \left( \begin{array}{cc} u_i(q) & 1 \\ v_i(q) & 0  \end{array} \right)
\end{eqnarray*}
and $u_i(q) := 1+q^{2i-1}$, $v_i(q) := (1+q^2)(1-q^{2(i-1)})/(1-q^2)$.

\begin{theorem}
For all $n\geq 2$,
\begin{eqnarray*}
  \TDIPD{n+1}{q} &=& \TDIPD{n}{q} + q^{2n} \TDIPO{n}{q} + \left(q^{2(n-1)}+\dfrac{q(1-q^{2n})}{1-q^2}\right) \TDIPD{n-1}{q} \\ 
  \TDIPO{n+1}{q} &=& \TDIPO{n}{q} + q^{2n} \TDIPD{n}{q} + \left(q^{2(n-1)}+\dfrac{q(1-q^{2n})}{1-q^2}\right) \TDIPO{n-1}{q} 
\end{eqnarray*}
with initial polynomials $\TDIPD{2}{q}, \TDIPO{2}{q} = 1+q+q^2$, $\TDIPD{3}{q} = (1+q+q^2+q^3)(1+q^3)+2q$ 
and $\TDIPO{3}{q} = (1+q+q^2+q^3)(1+q^3)+2q^5$.
\end{theorem}

\begin{proof}
Using Equation (\ref{dcon}),
\begin{eqnarray*}
  \TDIPD{n+2}{q} &=& \sum_{\pi \in I_{n+1}^D} q^{inv_D(\kodo{\pi}{n+2}{n+2})} 
  + \sum_{\pi \in I_{n+1}^O} q^{inv_D(\kodo{\pi}{-(n+2)}{n+2})} \\
  && + \sum_{k=1}^{n+1} \sum_{\pi \in I_n^D} q^{inv_D(\kodt{\pi}{k}{n+2})} + q^{inv_D(\kodt{\pi}{-k}{n+2})}.
\end{eqnarray*}
If $\pi \in I_{n+1}^D,I_{n+1}^O$, then $inv_D(\kodo{\pi}{n+2}{n+2}) = inv_D(\pi)$ 
and $inv_D(\kodo{\pi}{-(n+2)}{n+2}) = inv_D(\pi) +2(n+1)$.
Also if $\pi \in I_n^D$, then $inv_D(\kodt{\pi}{k}{n+2}) = 2n-2k+3+inv_D(\pi)$ and $inv_D(\kodt{\pi}{-k}{n+2})=inv_D(\pi)+2n$.
Hence,
\begin{eqnarray*}
  \TDIPD{n+2}{q} &=& \sum_{\pi \in I_{n+1}^D} q^{inv_D(\pi)} + \sum_{\pi \in I_{n+1}^O} q^{inv_D(\pi)+2(n+1)} \\
  && + \sum_{k=1}^{n+1} \sum_{\pi \in I_n^D} q^{inv_D(\pi)} (q^{2n-2k+3}+q^{2n}).
\end{eqnarray*}
The second recurrence is derived in the same manner by using Equation (\ref{ocon}).
\end{proof}

Let $J_n^D\subset I_n^D$ denote the class of all signed permutations such that $\pi(i) \neq \pm i$ for all $i \in [1,n]$ and consider
the generating function
$\mathcal{JD}_n(q):= \sum_{\pi \in J_n^D} q^{\inv_D (\pi)}$.

\begin{theorem}
\label{jdthm}
For all even $n\geq 2$, 
\begin{eqnarray*}
  \mathcal{JD}_{n}(q) &=& {2q^{n/2}} \prod_{i=1}^{n/2-1} \dfrac{(1+q^{4i})(1-q^{4i+2})}{1-q^2}.
\end{eqnarray*}
\end{theorem}

\begin{proof}
Since $J_n^D$ is a subclass of $I_n^D$ and from the characterization in Equation (\ref{dcon}), one has
\begin{eqnarray*}
  \mathcal{JD}_{n+4}(q) &=& \sum_{\pi \in J_{n+2}^D} q^{\inv_D (\moloko{\pi}{n+3}{n+4})}  + q^{\inv_D (\moloko{\pi}{-(n+3)}{n+4})}\\
  &&+ \sum_{1\leq i <j \leq n+2} \sum_{\pi \in J_n^D} \left(q^{\inv_D (\molokt{\pi}{i}{n+3}{j}{n+4})} 
    + q^{\inv_D (\molokt{\pi}{i}{n+4}{j}{n+3})} \right. \\
    && \qquad +q^{\inv_D (\molokt{\pi}{-i}{n+3}{j}{n+4})} + q^{\inv_D (\molokt{\pi}{-i}{n+4}{j}{n+3})} \\
    && \qquad +q^{\inv_D (\molokt{\pi}{i}{n+3}{-j}{n+4})} + q^{\inv_D (\molokt{\pi}{i}{n+4}{-j}{n+3})} \\
    && \left.\qquad + q^{\inv_D (\molokt{\pi}{-i}{n+3}{-j}{n+4})} + q^{\inv_D (\molokt{\pi}{-i}{n+4}{-j}{n+4})} \right).
\end{eqnarray*}
Now if $\pi \in J_{n+2}^D$ then $\inv_D (\moloko{\pi}{n+3}{n+4}) = \inv_D(\pi)+1$ and 
$\inv_D (\moloko{\pi}{-(n+3)}{n+4})= \inv_D(\pi)+4n+9$.
A careful analysis shows that for $\pi \in J_n^D$ and $1\leq i,j \leq n+2$, $i\neq j$, 
\begin{eqnarray*}
  \inv_D (\molokt{\pi}{i}{n+3}{j}{n+4}) &=& \inv_D(\pi)+4n-2(i+j) + 10+2\cdot\indicator{i>j}\\
  \inv_D (\molokt{\pi}{-i}{n+3}{-j}{n+4})&=& \inv_D(\pi)+4n+2(i+j) + 2-2\cdot \indicator{i>j}\\
  \inv_D (\molokt{\pi}{-i}{n+3}{j}{n+4}) &=& \inv_D(\pi) + 4n+2(i-j)+ 6 +2\cdot\indicator{i>j} \\
  \inv_D (\molokt{\pi}{i}{n+3}{-j}{n+4}) &=& \inv_D(\pi) + 4n+2(j-i)+ 6-2\cdot \indicator{i>j}.
\end{eqnarray*}
Thus we have
\begin{eqnarray*}
  \lefteqn{\mathcal{JD}_{n+4}(q)} \\ 
  &=& \mathcal{JD}_{n+2}(q) (q+q^{4n+9}) + \mathcal{JD}_{n}(q) \times \\ 
  && \sum_{1\leq i <j \leq n+2} \left( (q^{10}+q^{12})q^{4n-2(i+j)} + (q^{2}+1) q^{4n+2(i+j)} \right. \\
  && \qquad\qquad\qquad\left. +(q^6+q^8) q^{4n+2(i-j)} + (q^6+q^4) q^{4n+2(j-i)} \right)\\
  &=& \mathcal{JD}_{n+2}(q) (q+q^{4n+9}) + \mathcal{JD}_{n}(q) \dfrac{q^4(q^{4(n+1)}-1)(q^{2n+4}-1)(q^{2n}+1)}{(q^2-1)^2}.
\end{eqnarray*}
The result follows by inserting the expression from the theorem, we omit the details.
\end{proof}

Notice that if $n/2$ is even (resp. odd) then the coefficients of odd (resp. even) powers of $q$
in $\mathcal{JD}_{n}(q)$ are zero.

\begin{theorem}
The coefficients of the even (resp. odd) powers of $q$ in $\mathcal{JD}_{n}(q)$ are 
symmetric and unimodal when $n/2$ is even (resp. odd).
\end{theorem}

\begin{proof}
From Stanley~\cite[Proposition 1]{stanley.ann}, we have that if $A(q)$ and $B(q)$ are 
symmetric and unimodal polynomials, both with non-negative coefficients, then $A(q) B(q)$ is also symmetric and unimodal. 
From the expression in Theorem~\ref{jdthm}, one may write 
$\mathcal{JD}_{n+2}(q)=(q+q^2+\ldots +q^{2n}+ 2q^{2n+1}+ q^{2n+2}+\ldots + q^{4n+1})\mathcal{JD}_{n}(q)$.
The result follows inductively.
\end{proof}

The generating function of the descent polynomial over involutions of Coxeter groups of types $B_n$ and $D_n$ is
also seen to be symmetric, as was mentioned in Section 2, thanks to Hultman's~\cite{hultman} result.

\section{Comments}
Unlike the Eulerian polynomial, whose roots are all real and from which log-concavity of the coefficients follows,
the roots of all polynomials with the statistics mentioned above
are not real for $n\leq 14$.
Furthermore, they do not lie in the nice triangular $\pi /3$ region of the complex plane about the negative 
real-line from which it would be possible to infer log-concavity (see Stanley~\cite[Prop. 7]{stanley.ann}.)
Log-concavity of the coefficients holds numerically for all $n\leq 14$.
We extend the original conjecture,

\begin{conjecture}
\label{myconjecture}
For all $n\geq 4$, 
\renewcommand{\theenumi}{\roman{enumi}}\begin{enumerate}
\item the sequence $\{[q^i] {\mathcal{I}}_{n}^{maj}(q)\}_{i=0}^{n\choose 2}$ is log-concave,
\item for $2\leq i \leq {n\choose 2}-2$ (see Figure 1)
\begin{eqnarray*}
([q^i]\IP{n}{q})^2 &\geq& ([q^{i-2}]\IP{n}{q})([q^{i+2}] \IP{n}{q}),
\end{eqnarray*}
\item the sequences  $\{[q^{2i}] Inv_n^B(q)\}_{i\geq 0}$ and $\{[q^{2i+1}] Inv_n^B(q)\}_{i\geq 0}$ are unimodal,
\item the sequences  $\{[q^{2i}] Inv_n^D(q)\}_{i\geq 0}$ and $\{[q^{2i+1}] Inv_n^D(q)\}_{i\geq 0}$ are unimodal,
\item the sequences  $\{[q^{2i}] Inv_n^O(q)\}_{i\geq 0}$ and $\{[q^{2i+1}] Inv_n^O(q)\}_{i\geq 0}$ are unimodal.
\end{enumerate}\renewcommand{\theenumi}{\arabic{enumi}}
\end{conjecture}
We list here those polynomials for $n=10$ to exemplify these conjectures,
\begin{tiny}
\begin{eqnarray*}
\DP{10}{x} &=&  1 + 25 x + 289 x^{2} + 1397 x^{3} + 3036 x^{4} + 3036 x^{5} + 1397 x^{6} + 289 x^{7} + 25 x^{8} + x^{9}.\\
\EP{10}{x} &=&  1 + 45 x + 630 x^{2} + 3150 x^{3} + 4725 x^{4} + 945 x^{5}.\\
\MP{10}{x} &=&  1 + x + 2 x^{2} + 4 x^{3} + 7 x^{4} + 12 x^{5} + 19 x^{6} + 29 x^{7} + 44 x^{8} + 64 x^{9}  + 89 x^{10} + 119 x^{11}\\ 
&& + 158 x^{12} + 201 x^{13} + 250 x^{14} + 304 x^{15} + 358 x^{16} 
 + 412 x^{17} + 464 x^{18} + 508 x^{19} + 546 x^{20}\\ 
&&+572 x^{21}+584 x^{22} +584 x^{23}  + 572 x^{24} +546 x^{25} +508 x^{26} +464 x^{27} +412 x^{28} +358 x^{29}\\ 
&& + 304 x^{30}+250 x^{31}+201 x^{32}+158 x^{33}+119 x^{34}+89 x^{35}+64 x^{36}+44 x^{37}+29 x^{38}\\ 
&& + 19 x^{39}  + 12 x^{40} + 7 x^{41} + 4 x^{42} + 2 x^{43} + x^{44} + x^{45}.\\
\IP{10}{x} &=& 1+9 x+28 x^{2}+43 x^{3}+64 x^{4}+98 x^{5}+114 x^{6}+165 x^{7}+179 x^{8}+234 x^{9}+254 x^{10} \\ 
&& + 299 x^{11} + 333 x^{12} + 353 x^{13} + 408 x^{14} + 392 x^{15}  + 471 x^{16} + 411 x^{17} + 513 x^{18} + 409 x^{19}\\ 
&& + 529 x^{20} + 380 x^{21} + 517 x^{22}+335 x^{23} + 478 x^{24} + 281 x^{25} + 417 x^{26} + 225 x^{27} + 343 x^{28} \\ 
&& + 171 x^{29} + 264 x^{30} + 124 x^{31} + 189 x^{32} + 85 x^{33} + 123 x^{34} + 56 x^{35} + 72 x^{36}  + 35 x^{37} \\ 
&& + 37 x^{38} + 20 x^{39} + 16 x^{40} + 10 x^{41} + 5 x^{42} + 4 x^{43} + x^{44} + x^{45}.\\
\TDIPB{10}{x} &=&
1 +10\,x +36\,{x}^{2} +73\,{x}^{3} +157\,{x}^{4} +307\,{x}^{5} +456\,{x}^{6} +807\,{x}^{7} +1121\,{x}^{8} +1629\,{x}^{9}\\ 
&& +2323\,{x}^{10} +2835\,{x}^{11} +4124\,{x}^{12} +4508\,{x}^{13} +6468\,{x}^{14} +6715\,{x}^{15} +9256\,{x}^{16} \\ 
&& +9469\,{x}^{17} +12333\,{x}^{18} +12712\,{x}^{19} +15500\,{x}^{20} +16306\,{x}^{21} +18560\,{x}^{22} +20048\,{x}^{23}\\ 
\end{eqnarray*} \begin{eqnarray*}
&& +21334\,{x}^{24} +23730\,{x}^{25} +23626\,{x}^{26} +27127\,{x}^{27} +25285\,{x}^{28} +29989\,{x}^{29} +26242\,{x}^{30}\\ 
&& +32053\,{x}^{31} +26550\,{x}^{32} +33126\,{x}^{33} +26310\,{x}^{34} +33138\,{x}^{35} +25641\,{x}^{36} +32124\,{x}^{37} \\ 
&& +24639\,{x}^{38} +30194\,{x}^{39} +23393\,{x}^{40} +27534\,{x}^{41} +21953\,{x}^{42} +24364\,{x}^{43} +20369\,{x}^{44}\\ 
&& +20935\,{x}^{45}+18657\,{x}^{46}+17519\,{x}^{47}+16839\,{x}^{48}+14343\,{x}^{49}+14925\,{x}^{50}+11549\,{x}^{51}\\ 
&& +12956\,{x}^{52} +9199\,{x}^{53} +10967\,{x}^{54} +7288\,{x}^{55} +9019\,{x}^{56} +5762\,{x}^{57} +7178\,{x}^{58}\\ 
&& +4563\,{x}^{59} +5525\,{x}^{60} +3633\,{x}^{61} +4107\,{x}^{62} +2909\,{x}^{63} +2962\,{x}^{64} +2331\,{x}^{65}\\ 
&& +2084\,{x}^{66} +1858\,{x}^{67} +1444\,{x}^{68} +1460\,{x}^{69} +986\,{x}^{70} +1123\,{x}^{71} +671\,{x}^{72} +834\,{x}^{73}\\ 
&& +454\,{x}^{74} +589\,{x}^{75} +312\,{x}^{76} +394\,{x}^{77} +217\,{x}^{78} +255\,{x}^{79} +156\,{x}^{80} +156\,{x}^{81}\\ 
&& +111\,{x}^{82} +91\,{x}^{83} +79\,{x}^{84} +52\,{x}^{85} +56\,{x}^{86} +30\,{x}^{87} +40\,{x}^{88} +17\,{x}^{89} +26\,{x}^{90}\\ 
&& +10\,{x}^{91}+15\,{x}^{92}+5\,{x}^{93}+5\,{x}^{94}+3\,{x}^{95}+2\,{x}^{96}+3\,{x}^{97}+{x}^{98}+3\,{x}^{99}+{x}^{100}.\\
\TDIPD{10}{x} &=& 
1 +10x +35{x}^{2} +61{x}^{3} +97{x}^{4} +158{x}^{5} +204{x}^{6} +308{x}^{7} +370{x}^{8} +495{x}^{9} +595{x}^{10}\\ 
&& +734{x}^{11} +887{x}^{12} +1034{x}^{13} +1229{x}^{14} +1381{x}^{15} +1607{x}^{16} +1764{x}^{17} +2014{x}^{18} \\ 
&& +2182{x}^{19} +2432{x}^{20} +2606{x}^{21} +2827{x}^{22} +3012{x}^{23} +3175{x}^{24} +3377{x}^{25} +3451{x}^{26}\\ 
&& +3663{x}^{27} +3654{x}^{28} +3863{x}^{29} +3781{x}^{30} +3970{x}^{31} +3819{x}^{32} +3964{x}^{33} +3766{x}^{34}\\ 
&& +3859{x}^{35} +3642{x}^{36} +3670{x}^{37} +3432{x}^{38} +3402{x}^{39} +3156{x}^{40} +3085{x}^{41} +2844{x}^{42}\\ 
&& +2736{x}^{43} +2511{x}^{44} +2378{x}^{45} +2188{x}^{46} +2036{x}^{47} +1877{x}^{48} +1707{x}^{49} +1568{x}^{50}\\ 
&& +1396{x}^{51} +1284{x}^{52} +1128{x}^{53} +1035{x}^{54} +899{x}^{55} +818{x}^{56} +708{x}^{57} +642{x}^{58}\\ 
&& +553{x}^{59} +497{x}^{60} +428{x}^{61} +380{x}^{62} +322{x}^{63} +284{x}^{64} +236{x}^{65} +206{x}^{66} +168{x}^{67}\\ 
&& +142{x}^{68} +116{x}^{69} +98{x}^{70} +81{x}^{71} +68{x}^{72} +54{x}^{73} +46{x}^{74} +36{x}^{75} +32{x}^{76} +23{x}^{77}\\ 
&& +21{x}^{78} +18{x}^{79} +14{x}^{80} +11{x}^{81} +8{x}^{82} +5{x}^{83} +4{x}^{84} +2{x}^{85} +{x}^{86} +2{x}^{87} +{x}^{88}\\ 
&& +3{x}^{89} +{x}^{90}.\\
\TDIPO{10}{x} &=& 
1+8x+23x^2+41x^3+77x^4+120x^5 + 180x^6+268x^7+332x^8 +461x^9+547x^{10}\\ 
&& + 718x^{11}+835x^{12}+1040x^{13}+1181x^{14}+1407x^{15}+1569x^{16}+1808x^{17}+1994x^{18}\\ 
&& +2236x^{19}+2448x^{20}+2672x^{21}+2875x^{22}+3078x^{23}+3245x^{24}+3421x^{25} +3545x^{26}\\ 
&& +3679x^{27}+3758x^{28}+3847x^{29}+3877x^{30}+3926x^{31}+3899x^{32}+3906x^{33} +3826x^{34}\\ 
&&+3797x^{35}+3664x^{36}+3610x^{37}+3422x^{38}+3358x^{39}+3128x^{40}+3067x^{41}+2800x^{42}\\ 
&& +2744x^{43}+ 2461x^{44}+2408x^{45}+2138x^{46}+2080x^{47}+1835x^{48}+1759x^{49}+1542x^{50}\\ 
&& +1446x^{51}+1268x^{52}+1164x^{53}+1025x^{54}+919x^{55}+812x^{56}+708x^{57}+634x^{58}\\ 
&& +535x^{59}+491x^{60}+400x^{61}+374x^{62}+296x^{63} + 284x^{64}+218x^{65}+214x^{66}+156x^{67}\\ 
&& +150x^{68}+110x^{69} + 104x^{70}+81x^{71}+72x^{72}+54x^{73}+46x^{74}+34x^{75}+28x^{76}+23x^{77}\\ 
&& +17x^{78}+16x^{79}+10x^{80}+11x^{81}+6x^{82}+7x^{83}+4x^{84}+4x^{85}+x^{86}+2x^{87}+x^{88}\\ 
&& +x^{89}+x^{90}.
\end{eqnarray*}
\end{tiny}
\begin{figure}
\thediagram
\caption{The coefficients $[x^{2i}]\IP{10}{x}$, $[x^{2i+1}]\IP{10}{x}$}
\end{figure}

\section*{Acknowledgments}
The author would like to thank M. Bousquet-M\'elou and 
F. Brenti for helpful comments
and also D. Foata and V. Strehl for pointers to the relevant literature.

\end{document}